\documentclass[12pt]{article}
\usepackage{amsfonts}
\usepackage{amsmath}

\newtheorem{theorem}{Theorem}

\newtheorem{corollary}[theorem]{Corollary}

\newtheorem{definition}[theorem]{Definition}
\newtheorem{example}[theorem]{Example}

\newtheorem{problem}[theorem]{Problem}
\newtheorem{proposition}[theorem]{Proposition}

\newenvironment{proof}[1][Proof]{\noindent\textbf{#1.} }{\ \rule{0.5em}{0.5em}}
\input{tcilatex}

\begin{document}

\title{Partial Transpose of Permutation Matrices}
\author{Qing-Hu Hou\thanks{\texttt{hou@nankai.edu.cn}} \\
Center for Combinatorics, Nankai University, \\
Tianjin 300071, P. R. China \and Toufik Mansour\thanks{\texttt{%
toufik@math.haifa.ac.il}} \\
Department of Mathematics, University of Haifa, \\
Haifa 31905, Israel \and Simone Severini\thanks{\texttt{simoseve@gmail.com}} 
\\
Institute for Quantum Computing and\\
Department of Combinatorics \& Optimization, \\
University of Waterloo, \\
Waterloo N2L 3G1, Canada}
\maketitle

\begin{abstract}
The partial transpose of a block matrix $M$ is the matrix obtained by
transposing the blocks of $M$ independently. We approach the notion of
partial transpose from a combinatorial point of view. In this perspective,
we solve some basic enumeration problems concerning the partial transpose of
permutation matrices. More specifically, we count the number of permutations
matrices which are equal to their partial transpose and the number of
permutation matrices whose partial transpose is still a permutation. We
solve these problems also when restricted to symmetric permutation matrices
only.
\end{abstract}

\section{Introduction}

The \emph{partial transpose} (or, equivalently, \emph{partial transposition}%
) is a linear algebraic concept, which can be interpreted as a simple
generalization of the usual matrix transpose. In the present paper, we
consider partial transpose from a combinatorial point of view. More
specifically, we solve some enumeration problems concerning the partial
transpose of permutation matrices.

Even if this notion is a natural one, to the knowledge of the authors, it
has never been directly studied by the linear algebra community. On the
other hand, partial transpose is an important tool in the mathematical
theory of quantum entanglement. For this reason, partial transpose appears
often in works contextual with quantum information theory. We will spend a
few paragraphs on this, just for taking a snapshot of the scenario in which
this notion arises.

Bru\ss\ and Macchiavello \cite{bm} give an excellent explanation of the
meaning of partial transpose in quantum information theory. Its primary use
is materialized in the so-called \emph{PPT-criterion}, where
\textquotedblleft PPT\textquotedblright\ stands for \emph{Positive Partial
Transpose}. The criterion, firstly discovered by Peres \cite{p} and the
Horodeckis \cite{h} (see also \cite{nc}), is as follows: if the density
matrix (or, equivalently, the state) of a quantum mechanical system with
composite dimension $pq$ is entangled, with respect to the subsystems of
dimension $p$ and $q$, then its partial transpose is positive. The converse
of the implication is not necessarily true. However, under certain
restrictions, for example, when the dimension of the density matrix is six,
the PPT-criterion is necessary and sufficient.

There is a number of problems suggested by the PPT-criterion. In particular,
in order to shed light onto the structure of the set of density matrices, it
would be important to characterize those for which the criterion is valid.
An open question of practical importance is to prove or disprove that
certain states, which are said to be non-distillable, have positive partial
transpose. However, there is strong evidence that there exist
non-distillable states with negative partial transpose, which would be then
called NPT-bound entangled states. Regarding this topic, see the important
references \cite{di, du}, or \cite{ch}, for an account on recent discussions.

Looking at the notion of partial transpose from the combinatorial point of
view is an appealing topic, because it has the potential to uncover patters
in the set of density matrices and indicate connections with other
mathematical objects, and this may turn out to be helpful in understanding
physical properties. As a matter of fact there have been a number of recent
papers considering entanglement in discrete settings (see, \emph{e.g.}, \cite%
{b, g, k}).

Here we state and solve some basic enumeration problems involving partial
transpose of permutation matrices. Permutations appear in fact to be a
simple, yet a rich territory to explore. Enumeration is a good first step
towards the quantitative understanding of the structure of a set.

In particular, we count the number of permutations matrices which are equal
to their partial transpose and the number of permutation matrices whose
partial transpose is still a permutation. We solve these problems also when
restricted to symmetric permutation matrices only (\emph{i.e.}, induced by
involutions).

Apart from considerations related to symmetry, given that symmetry often
predisposes to relations between different combinatorial objects, a further
reason to look at involutions comes from \cite{b}. A permutation matrix
associated to a involution without fixed points can be seen as the adjacency
matrix of the disjoint union of matchings and self-loops. Since the
combinatorial Laplacian of any graph is a density matrix after appropriate
normalization \cite{b}, counting the number of involutions whose partial
transpose is a permutation, is equivalent to count the number of these
states with positive partial transpose. However, the PPT-criterion is not
sufficient also for this extremely restricted class. There actually are
disconnected graphs whose Laplacian is entangled even if its partial
transpose is positive \cite{hi}.

\bigskip

The organization of the paper is as follows. In the next section, we give
the required definitions and formally state our problems. Section \ref{sec3}
deals with permutations whose partial transpose is a permutation; Section %
\ref{sec4}, with permutations equal to their partial transpose; Section \ref%
{sec5}, with involutions whose partial transpose is a permutation

\section{Definitions, statements of the problems and examples}

The following is a formal definition of the partial transpose of a matrix:

\begin{definition}
Let $M$ be an $n\times n$ matrix with real entries. Let us assume that $n=pq$%
, where $p$ and $q$ are chosen arbitrarily. Under this assumption, we can
look at the matrix $M$ as partitioned into $p^{2}$ blocks each one $q\times
q $. The \emph{partial transpose} of $M$, denoted by $M^{\Gamma _{p}}$, is
the matrix obtained from $M$, by transposing independently each of its $%
p^{2} $ blocks. Formally, if%
\begin{equation*}
M=\left( 
\begin{array}{ccc}
\mathcal{B}_{1,1} & \cdots & \mathcal{B}_{1,p} \\ 
\vdots & \ddots & \vdots \\ 
\mathcal{B}_{p,1} & \cdots & \mathcal{B}_{p,p}%
\end{array}%
\right)
\end{equation*}%
then%
\begin{equation*}
M^{\Gamma _{p}}=\left( 
\begin{array}{ccc}
\mathcal{B}_{1,1}^{T} & \cdots & \mathcal{B}_{1,p}^{T} \\ 
\vdots & \ddots & \vdots \\ 
\mathcal{B}_{p,1}^{T} & \cdots & \mathcal{B}_{p,p}^{T}%
\end{array}%
\right) ,
\end{equation*}%
where $\mathcal{B}_{i,j}^{T}$ denotes the transpose of the block $\mathcal{B}%
_{i,j}$, for $1\leq i,j\leq p$.
\end{definition}

Notice that, by taking the adjoint $\mathcal{B}_{i,j}^{\dagger }$, instead
of the transpose $\mathcal{B}_{i,j}^{T}$, the notion of partial transpose
can be easily extended to matrices with complex entries. This is something
which we will not need here. Note that we have defined partial transpose
with respect to the parameter $p$. We could have also defined partial
transpose with respect to the parameter $q$, by treating the blocks of $M$
as the entries of a $p\times p$ matrix. Formally, 
\begin{equation*}
M^{\Gamma _{q}}=\left( 
\begin{array}{ccc}
\mathcal{B}_{1,1} & \cdots & \mathcal{B}_{p,1} \\ 
\vdots & \ddots & \vdots \\ 
\mathcal{B}_{1,p} & \cdots & \mathcal{B}_{p,p}%
\end{array}%
\right) .
\end{equation*}%
That is, the block $\mathcal{B}_{i,j}$ in $M$ is the block $\mathcal{B}%
_{j,i} $ in $M^{\Gamma _{q}}$, for all $1\leq i,j\leq p$. The term
\textquotedblleft partial transpose\textquotedblright\ also indicates the
actual operation required to obtain the matrix partial transpose as defined
here.

We will consider partial transpose of permutation matrices. Let us recall
that a \emph{permutation matrix} of \emph{size} $n$ is an $n\times n$
matrix, with entries in the set $\{0,1\},$ such that each row and each
column contains exactly one nonzero entry. A \emph{permutation} of \emph{%
length} $n$ is a bijection $\pi :[n]\longrightarrow \lbrack n]$, where $%
[n]=\{1,2,...,n\}$. Given an $n\times n$ permutation matrix $P$, there is a
unique permutation $\pi $ of length $n$ associated to $P$, such that $\pi
(i)=j$ if and only if $P_{i,j}=1$. Let us denote by $S_{n}$ the set of all $%
n\times n$ permutation matrices. With an innocuous abuse of notation, we
write $S_{n}$ also for the set of all permutations of length $n$.

In standard linear notation, a permutation $\pi \in S_{n}$ can be written as
a word of the form $\pi (1)\pi (2)...\pi (n)$. It may be interesting to
point out that a permutation and its partial transpose share a common
property, related to the sum of the row indices. The cells in the table
below contain ordered pairs: each left element of the pairs is a permutation 
$\pi \in S_{4}$; each right element is the ordered list of the row indices
of the one entries in the matrix $P^{\Gamma _{2}}$:

\begin{equation*}
\begin{tabular}{||c|c|c|c|c|c||}
\hline\hline
$1234,1234$ & $1243,1243$ & $1324,1414$ & $1342,1432$ & $1423,1441$ & $%
1432,1432$ \\ \hline
$2134,2134$ & $2143,2143$ & $2314,4114$ & $2341,4123$ & $2413,1414$ & $%
2431,1423$ \\ \hline
$3142,2314$ & $3142,2323$ & $3214,3214$ & $3241,3223$ & $3412,3412$ & $%
3421,3421$ \\ \hline
$4123,2341$ & $4132,2332$ & $4213,2314$ & $4231,2323$ & $4312,4312$ & $%
4321,4321$ \\ \hline\hline
\end{tabular}%
\end{equation*}

\begin{example}
\emph{If}%
\begin{equation*}
P=\left( 
\begin{array}{cccc}
0 & 0 & 1 & 0 \\ 
1 & 0 & 0 & 0 \\ 
0 & 0 & 0 & 1 \\ 
0 & 1 & 0 & 0%
\end{array}%
\right) ,
\end{equation*}%
\emph{then }%
\begin{equation*}
P^{\Gamma _{2}}=\left( 
\begin{array}{cccc}
0 & 1 & 1 & 0 \\ 
0 & 0 & 0 & 0 \\ 
0 & 0 & 0 & 0 \\ 
0 & 1 & 1 & 0%
\end{array}%
\right) .
\end{equation*}%
\emph{The matrix }$P$\emph{\ is induced by the permutation }$\pi =3124$\emph{%
. The ordered list of the row indices of the one entries in }$P^{\Gamma
_{2}} $\emph{\ is }$2323$\emph{.}
\end{example}

For every permutation $\pi \in S_{n}$, where $n=pq$, we have 
\begin{eqnarray}
\sum_{i=1}^{n}\pi (i)
&=&\sum_{i=1}^{n}i=\sum_{P_{i,j}=1}i=\sum_{P_{i,j}^{\Gamma _{p}}=1}i
\label{gau} \\
&=&n(n+1)/2.  \notag
\end{eqnarray}%
This is straightforward. Let the $(ap+i,b(a,i)p+j(a,i))$-th entry of $P$ be
equal to $1$. Then $b(a,i)$ runs $p$ times over $0,\ldots ,q-1$ and $j(a,i)$
runs $q$ times over $1,\ldots ,p$. Thus, 
\begin{equation*}
\sum_{a,i}b(a,i)=p\binom{q}{2}
\end{equation*}%
and 
\begin{equation*}
\sum_{a,i}j(a,i)=q\binom{p+1}{2}.
\end{equation*}%
Therefore 
\begin{equation*}
\sum_{a,i}ap+j(a,i)=p\binom{q}{2}+q\binom{p+1}{2}=\binom{n+1}{2},
\end{equation*}%
which validates Eq. (\ref{gau}).

Let us recall that a permutation matrix $P$ is said to be a \emph{involution}
if $P=P^{T}$ and $P$ is not the identity matrix. We will solve the following
problems:

\begin{problem}
\label{pro2}Count the number of permutation matrices $P\in S_{pq}$ such that 
$P^{\Gamma _{p}}\in S_{pq}$.
\end{problem}

\begin{example}
\emph{When }$p=q=2$\emph{, we have all together }$12$\emph{\ matrices }$P\in
S_{4}$\emph{\ such that }$P^{\Gamma _{2}}\in S_{4}$\emph{. Among these, }$8$%
\emph{\ are the block-matrices of the forms }%
\begin{equation}
\begin{tabular}{lll}
$\left( 
\begin{array}{cc}
\ast  & \mathbf{0} \\ 
\mathbf{0} & \ast 
\end{array}%
\right) $ & and & $\left( 
\begin{array}{cc}
\mathbf{0} & \ast  \\ 
\ast  & \mathbf{0}%
\end{array}%
\right) $%
\end{tabular}%
.  \label{litt}
\end{equation}%
\emph{The remaining }$4$\emph{\ matrices are}%
\begin{equation}
\left( 
\begin{array}{cccc}
1 & 0 & 0 & 0 \\ 
0 & 0 & 0 & 1 \\ 
0 & 0 & 1 & 0 \\ 
0 & 1 & 0 & 0%
\end{array}%
\right) ,\left( 
\begin{array}{cccc}
0 & 1 & 0 & 0 \\ 
0 & 0 & 1 & 0 \\ 
0 & 0 & 0 & 1 \\ 
1 & 0 & 0 & 0%
\end{array}%
\right) ,\left( 
\begin{array}{cccc}
0 & 0 & 1 & 0 \\ 
0 & 1 & 0 & 0 \\ 
1 & 0 & 0 & 0 \\ 
0 & 0 & 0 & 1%
\end{array}%
\right) ,\left( 
\begin{array}{cccc}
0 & 0 & 0 & 1 \\ 
1 & 0 & 0 & 0 \\ 
0 & 1 & 0 & 0 \\ 
0 & 0 & 1 & 0%
\end{array}%
\right) .  \label{mat}
\end{equation}
\end{example}

\begin{problem}
\label{pro3}Count the number of permutation matrices $P\in S_{pq}$ such that 
$P^{\Gamma _{p}}=P$.
\end{problem}

\begin{example}
\emph{When }$p=q=2$\emph{, we have all together }$10$\emph{\ matrices }$P\in
S_{4}$\emph{\ such that }$P^{\Gamma _{2}}=P$\emph{. Among these, }$8$\emph{\
are the block matrices in Eq. (\ref{litt}). The remaining }$2$\emph{\
matrices are the first and the third matrix in Eq. (\ref{mat}).}
\end{example}

\begin{problem}
\label{pro4}Count the number of involutions $P\in S_{pq}$ such that $%
P^{\Gamma _{p}}\in S_{pq}$.
\end{problem}

\begin{example}
\emph{When }$p=q=2$\emph{, we have all together }$8$\emph{\ involutions }$%
P\in S_{4}$\emph{\ such that }$P^{\Gamma _{2}}=P$\emph{: }%
\begin{eqnarray*}
&&\left( 
\begin{array}{cccc}
1 & 0 & 0 & 0 \\ 
0 & 0 & 0 & 1 \\ 
0 & 0 & 1 & 0 \\ 
0 & 1 & 0 & 0%
\end{array}%
\right) ,\left( 
\begin{array}{cccc}
0 & 0 & 1 & 0 \\ 
0 & 1 & 0 & 0 \\ 
1 & 0 & 0 & 0 \\ 
0 & 0 & 0 & 1%
\end{array}%
\right) ,\left( 
\begin{array}{cccc}
0 & 1 & 0 & 0 \\ 
1 & 0 & 0 & 0 \\ 
0 & 0 & 0 & 1 \\ 
0 & 0 & 1 & 0%
\end{array}%
\right) ,\left( 
\begin{array}{cccc}
0 & 0 & 1 & 0 \\ 
0 & 0 & 0 & 1 \\ 
1 & 0 & 0 & 0 \\ 
0 & 1 & 0 & 0%
\end{array}%
\right) \\
&&\left( 
\begin{array}{cccc}
1 & 0 & 0 & 0 \\ 
0 & 1 & 0 & 0 \\ 
0 & 0 & 1 & 0 \\ 
0 & 0 & 0 & 1%
\end{array}%
\right) ,\left( 
\begin{array}{cccc}
0 & 0 & 0 & 1 \\ 
0 & 0 & 1 & 0 \\ 
0 & 1 & 0 & 0 \\ 
1 & 0 & 0 & 0%
\end{array}%
\right) ,\left( 
\begin{array}{cccc}
1 & 0 & 0 & 0 \\ 
0 & 1 & 0 & 0 \\ 
0 & 0 & 0 & 1 \\ 
0 & 0 & 1 & 0%
\end{array}%
\right) ,\left( 
\begin{array}{cccc}
0 & 1 & 0 & 0 \\ 
1 & 0 & 0 & 0 \\ 
0 & 0 & 1 & 0 \\ 
0 & 0 & 0 & 1%
\end{array}%
\right) .
\end{eqnarray*}
\end{example}

\section{Permutations whose partial transpose is a permutation\label{sec3}}

We begin with Problem \ref{pro2}. For a permutation matrix $P\in S_{pq}$,
let us denote by $\mathcal{B}_{i,j}$ the block located in the $i$-th row and 
$j$-th column. Let further $A_{i,j},B_{i,j}\subseteq \lbrack q]=\{1,2,\ldots
,q\}$ be the sets of relative row indices and column indices of the $1$'s in
the block $\mathcal{B}_{i,j}$. For example, given 
\begin{equation*}
\left( 
\begin{array}{cccc}
0 & 0 & 0 & 1 \\ 
1 & 0 & 0 & 0 \\ 
0 & 1 & 0 & 0 \\ 
0 & 0 & 1 & 0%
\end{array}%
\right) ,
\end{equation*}%
we have 
\begin{equation*}
A_{1,1}=\{2\},A_{1,2}=\{1\},A_{2,1}=\{1\},A_{2,2}=\{2\},
\end{equation*}%
and 
\begin{equation*}
B_{1,1}=\{1\},B_{1,2}=\{2\},B_{2,1}=\{2\},B_{2,2}=\{1\}.
\end{equation*}%
Clearly, $A_{i,j}$ has the same cardinality as $B_{i,j}$, which we denote by 
$r_{i,j}$. For fixed $A_{i,j}$ and $B_{i,j}$, we have $r_{i,j}!$ ways to
place $1$'s in $\mathcal{B}_{i,j}$. Therefore, the number of required
matrices equals the number of $A_{i,j},B_{i,j}$'s multiplied by $%
\prod_{i,j}r_{i,j}!$. At this stage, we impose the required constraints on $%
A_{i,j}$ and $B_{i,j}$. We know that $P$ is a permutation matrix if and only
if 
\begin{equation}
A_{i,j}\cap A_{i,k}=\emptyset ,\quad \text{for every}\,i,j,k\text{ with }%
j\not=k,  \label{Aij}
\end{equation}%
\begin{equation}
B_{i,j}\cap B_{k,j}=\emptyset ,\quad \text{for every}\,i,j,k\text{ with }%
i\not=k,  \label{Bji}
\end{equation}%
and 
\begin{align}
\bigcup_{j=1}^{p}A_{i,j}& =[q],\quad \text{for }i=1,2,...,p,  \label{Uij} \\
\bigcup_{i=1}^{p}B_{i,j}& =[q],\quad \text{for }j=1,2,...,p.
\end{align}%
We need that $P^{\Gamma _{p}}$ is also a permutation matrix. Therefore, we
have 
\begin{equation}
A_{i,j}\cap A_{k,j}=\emptyset ,\quad \text{for every }i,j,k\text{ with }%
i\not=k,  \label{Aji}
\end{equation}%
\begin{equation}
B_{i,j}\cap B_{i,k}=\emptyset ,\quad \text{for every}\,i,j,k\text{ with }%
j\not=k,  \label{Bij}
\end{equation}%
and 
\begin{align}
\bigcup_{i=1}^{p}A_{i,j}& =[q],\quad \text{for }j=1,2,...,p, \\
\bigcup_{j=1}^{p}B_{i,j}& =[q],\quad \text{for }i=1,2,...,p.  \label{Uji}
\end{align}%
Let 
\begin{equation*}
A_{\pi }=\bigcap_{i=1}^{p}A_{i,\pi _{i}}\quad \text{and}\quad B_{\pi
}=\bigcap_{i=1}^{p}B_{i,\pi _{i}},\quad \forall \,\pi \in S_{p},
\end{equation*}%
By Eqs. \eqref{Aij}--\eqref{Uij}, we know that 
\begin{equation}
A_{i,j}=\bigcup_{\pi _{i}=j}A_{\pi },\quad B_{i,j}=\bigcup_{\pi
_{i}=j}B_{\pi }.  \label{AB-pi}
\end{equation}%
From \eqref{Aij} and \eqref{Bji}, we can then write 
\begin{equation}
A_{\pi }\cap A_{\sigma }=B_{\pi }\cap B_{\sigma }=\emptyset ,\quad \text{for
every }\pi ,\sigma \in S_{p}\text{ with }\pi \not=\sigma .  \label{cap}
\end{equation}%
Furthermore, 
\begin{equation}
\bigcup_{\pi \in S_{p}}A_{\pi }=\bigcup_{\pi \in S_{p}}B_{\pi }=[q].
\label{cup}
\end{equation}%
Conversely, given two set partitions $\{A_{\pi }\}$ and $\{B_{\pi }\}$ of $%
[q]$, satisfying Eqs. \eqref{cap} and \eqref{cup}, we may define $A_{i,j}$
and $B_{i,j}$ by Eq. \eqref{AB-pi}. One can easily check that Eqs. %
\eqref{Aij}--\eqref{Uji} hold. The only restriction on the $A_{\pi }$'s and
the $B_{\pi }$'s is that the cardinalities of $A_{i,j}$ and $B_{i,j}$ should
be the equal. Let $a_{\pi }$ and $b_{\pi }$ denote the cardinalities of $%
A_{\pi }$ and $B_{\pi }$, respectively. On the basis of the above lines, we
can state the following result:

\begin{theorem}
\label{ppa}Let $Z(p,q)$ the number of permutation matrices $P\in S_{pq}$
such that $P^{\Gamma _{p}}\in S_{pq}$. Then 
\begin{equation}
Z(p,q)=\sum_{\substack{ \sum a_{\pi }=\sum b_{\pi }=q  \\ \sum_{\pi
_{i}=j}a_{\pi }=\sum_{\pi _{i}=j}b_{\pi }}}\frac{q!^{2}}{\prod_{\pi }a_{\pi
}!b_{\pi }!}\prod_{i,j=1}^{p}\left( \sum_{\pi _{i}=j}a_{\pi }\right) !,
\label{f-1}
\end{equation}%
where the sum runs over all $a_{\pi },b_{\pi }\in \mathbb{Z}$.
\end{theorem}

The following corollary shows a neat expression for the special case $P\in
S_{2q}$:

\begin{corollary}
The number of permutation matrices $P\in S_{2q}$ such that $P^{\Gamma
_{2}}\in S_{2q}$ is 
\begin{equation*}
Z(2,q)=q!(q+1)!.
\end{equation*}
\end{corollary}

The pattern avoidance language is now a standard tool for characterizing
classes of permutations (see \cite{we}). It would be natural to find a
characterization of the set of permutations given in Theorem \ref{ppa} in
terms of pattern avoidance.

\section{Permutations equal to their partial transpose\label{sec4}}

We now focus on Problem \ref{pro3}. We then ask that $\mathcal{B}_{i,j}=%
\mathcal{B}_{i,j}^{T}$. Hence, $A_{i,j}=B_{i,j}$. Additionally, given $%
A_{i,j}$, the number of ways to put $1$'s in the block $\mathcal{B}_{i,j}$
is exactly the number of involutions of length $q$, which we denote by $I(q)$%
. It is well-known that (see, \emph{e.g.}, \cite{s}, Example 5.2.10)%
\begin{equation}
I(q)=\sum_{\substack{ j=0  \\ j\text{ even}}}^{q}\binom{q}{j}\frac{j!}{%
2^{j/2}(j/2)!}  \label{invnum}
\end{equation}%
and 
\begin{equation*}
I(q+1)=I(q)+q\cdot I(q-1).
\end{equation*}%
With the same analysis carried on for Theorem \ref{ppa}, we can directly
obtain the number of desired matrices:

\begin{theorem}
\label{ppe}Let $Z_{e}(p,q)$ the number of permutation matrices $P\in S_{pq}$
such that $P=P^{\Gamma _{p}}$. Then%
\begin{equation}
Z_{e}(p,q)=\sum_{\sum a_{\pi }=q}\frac{q!}{\prod_{\pi }a_{\pi }!}%
\prod_{i,j=1}^{p}i\left( \sum_{\pi _{i}=j}a_{\pi }\right) ,
\end{equation}%
where the sum runs over all $a_{\pi }\in \mathbb{Z}$, with $\pi \in S_{p}$.
\end{theorem}

When taking $p=2$, the number of permutation matrices is given in the next
corollary:

\begin{corollary}
The number of permutation matrices $P\in S_{2q}$ such that $P=P^{\Gamma
_{2}} $ is%
\begin{equation*}
Z_{e}(2,q)=\sum_{r=0}^{q}\binom{q}{r}^{2}I(r)^{2}I(q-r)^{2}.
\end{equation*}
\end{corollary}

\section{Involutions whose partial transpose is a permutation\label{sec5}}

In this section, we present a solution of Problem \ref{pro4}. Let $P$ be the
involution defined by the ordered pairs $(aq+i,bq+j)$, where $0\leq a,b\leq
p-1$, $1\leq i,j\leq q$ and $(a,i)\not=(b,j)$. Note that the partial
transpose keeps fixed the $1$'s in the diagonal. So, the only possible
permutation matrices after partial transpose would be the identity matrix $%
\mathrm{Id}$ or $P$ itself. In the first case, we must have $P=\mathrm{Id}$,
since we get back to the original matrix by applying twice the partial
transpose operation. Therefore, we only need to consider the second case,
that is, when $P$ remains invariant under partial transpose. Notice that the 
$(aq+i,bq+j)$-th and the $(bq+j,aq+i)$-th entry of the permutation matrix
are $1$'s. After partial transpose, the $(aq+j,bq+i)$-th and the $%
(bq+i,aq+j) $-th entry are $1$'s. Thus we have 
\begin{eqnarray*}
(aq+i,bq+j) &=&(aq+j,bq+i), \\
(bq+j,aq+i) &=&(bq+i,aq+j),
\end{eqnarray*}%
or 
\begin{eqnarray*}
(aq+i,bq+j) &=&(bq+i,aq+j), \\
(bq+j,aq+i) &=&(aq+j,bq+i).
\end{eqnarray*}%
That is, $i=j$ or $a=b$. Hence, the desired involutions are of type $%
(aq+i,aq+j)$, with $i\not=j$, or, of type $(aq+i,bq+i)$, with $a\not=b$. We
can then state the following fact:

\begin{theorem}
Let $Z_{t}(p,q)$ be the number of involutions $P\in S_{pq}$ such that $%
P^{\Gamma _{p}}\in S_{pq}$, or, equivalently, $P^{\Gamma _{p}}=P$. Then 
\begin{equation*}
Z_{t}(p,q)=2p\binom{q}{2}+2q\binom{p}{2}.
\end{equation*}
\end{theorem}

\begin{corollary}
\label{tre}The following statements hold true:

\begin{itemize}
\item $Z_{t}(q+1,q)=q(q+1)(2q-1)$;

\item $Z_{t}(q,q)=2(q^{3}-q^{2})$.
\end{itemize}
\end{corollary}

The numbers $Z_{t}(q+1,q)/2$ are called \emph{octagonal pyramidal numbers},
and count the ways of covering a $2q\times 2q$ lattice with $2q^{2}$
dominoes with exactly $2$ horizontal dominoes (\cite{sl}, Seq. A002414). The
numbers $Z_{t}(q,q)$ count the possible rook moves on an a $q\times q$
chessboard (\cite{b, sl}, Seq. A002414).

To conclude this section, even if these are simple facts, it may be
clarifying to remark the following:

\begin{proposition}
The following statements hold true for all $p$ and $q$:

\begin{itemize}
\item $Z(p,q)=Z(q,p)$;

\item In general, $Z_{e}(p,q)\not=Z_{e}(p,q)$;

\item $Z_{t}(p,q)=Z_{t}(q,p)$.
\end{itemize}
\end{proposition}

\begin{proof}
While the second point is obvious, the other two can be verified by the
following bijection. Suppose that the $(ap+i,b(a,i)p+j(a,i))$-th entry of $P$
is $1$. Then let the $((i-1)q+(a+1),(j(a,i)-1)q+(b(a,i)+1))$-th entry of $%
P^{\prime }$ be $1$. If the partial transpose of $P$ is a permutation, then $%
ap+j(a,i)$ and $b(a,i)p+i$ run from $1$ to $n$, for $0\leq a\leq q-1,1\leq
i\leq p$. Thus, $(i-1)q+(b(a,i)+1)$ and $j(a,i)-1)q+(a+1)$ run from $1$ to $%
n $ also. This implies that the partial transpose of $P^{\prime }$ is a
permutation.
\end{proof}

\bigskip

\noindent \emph{Acknowledgment}. This paper arised from discussions, while
the third author was visiting the Center for Combinatorics at Nankai
University. The financial support of this institution is gratefully
acknowledged.

\end{document}